\documentclass{svjour2}
\smartqed
\usepackage{amsmath,amssymb}
\usepackage{graphicx}
\usepackage{booktabs}
\usepackage{tikz}
\usepackage{url}
\usepackage[hidelinks]{hyperref}

\journalname{Journal of Statistical Physics}
\makeatletter\def\subclassname{{\bfseries Mathematics Subject Classification (2020)}\enspace}\makeatother

\newcommand{\Tr}{\operatorname{Tr}}
\newcommand{\A}{\mathcal A}
\newcommand{\N}{\mathbb N}
\newcommand{\Q}{\mathbb Q}
\newcommand{\C}{\mathcal C}
\newcommand{\Pc}{\mathcal P}

\begin{document}

\title{On the lower bound of the refined Bessis--Moussa--Villani conjecture: proof of the cases $S_{6,3}$ and $S_{7,3}$ for $3\times3$ matrices}

\author{Daniel H\"agele}

\institute{D. H\"agele \at
Experimental Physics VI (AG), Faculty of Physics and Astronomy, Ruhr University Bochum, Universit\"atsstra\ss e 150, D-44780 Bochum, Germany\\
\email{daniel.haegele@rub.de}}

\date{Received: date / Accepted: date}

\maketitle

\begin{abstract}
For Hermitian matrices $\alpha,\beta$ let $\A_{n,m}(e^\alpha,e^\beta)$ be the average of the traces of all words with $n$ letters $e^\alpha$ and $m$ letters $e^\beta$. The refined Bessis--Moussa--Villani conjecture contains the lower bound $\A_{n,m}(e^\alpha,e^\beta)\ge\Tr e^{n\alpha+m\beta}$ (Tr denotes the trace), which for $n=m=1$ is the Golden--Thompson inequality. We prove this lower bound for all $n+m\le7$ when the matrices have size at most three, and for $\min(n,m)\le2$ in every dimension. In all these cases we prove the stronger inequality $\A_{n,m}(A,B)\ge\Tr(A^{n/m}B)^m$ for positive semidefinite $A,B$, or its mirror image with the roles of $A$ and $B$ exchanged; the lower bound of the refined conjecture follows from the Araki--Lieb--Thirring inequality. The cases with three letters $B$ in dimension three, $(n,m)=(3,3)$ and $(4,3)$, are reduced by an exact identity and the inequality of arithmetic and geometric means to a polynomial inequality in two variables, which is certified in exact rational arithmetic by Bernstein expansion on three boxes. The certificate is a short program in the ancillary files.
\end{abstract}

\section{Introduction}\label{sec:intro}
In 1975 Bessis, Moussa and Villani conjectured that $t\mapsto\Tr e^{H-tV}$ ($\Tr$ denotes the trace) is the Laplace transform of a positive measure for Hermitian $H$ and positive semidefinite $V$ \cite{bmv1975}. The motivation came from quantum statistical mechanics: the conjecture yields monotone Pad\'e bounds for partition functions. Lieb and Seiringer showed that it is equivalent to the statement that for positive semidefinite $A,B$ all coefficients of the polynomial $t\mapsto\Tr(A+tB)^p$ are nonnegative \cite{liebseiringer2004}. The coefficient of $t^m$ is $\Tr S_{p,m}(A,B)$, where $S_{p,m}(A,B)$ denotes the sum of all $\binom pm$ words of length $p$ with $m$ letters $B$ and $p-m$ letters $A$. After partial results, among them the cases $p\le7$ \cite{haegele2007} and $p\le13$ \cite{klepschweighofer2008}, the conjecture was proved by Stahl \cite{stahl2013}.

For commuting matrices all words in $S_{p,m}$ coincide. The author proposed a quantitative refinement, listed as Problem~40 of the Open Quantum Problems collection \cite{oqp40}. With
\[
\A_{n,m}(A,B):=\binom{n+m}{m}^{-1}\Tr S_{n+m,m}(A,B)
\]
it reads, for Hermitian $\alpha,\beta$ of any size $d$,
\begin{equation}\label{eq:refined}
\Tr e^{n\alpha+m\beta}\ \le\ \A_{n,m}(e^\alpha,e^\beta)\ \le\ \Tr\big(e^{n\alpha}e^{m\beta}\big).
\end{equation}
For $n=m=1$ the left inequality is the Golden--Thompson inequality \cite{golden1965,thompson1965}, so \eqref{eq:refined} interpolates between Golden--Thompson and Stahl's theorem. The upper bound in \eqref{eq:refined} is false: Cha and Lee found $3\times3$ counterexamples for $n=m=5$ \cite{cha2026}; see also \cite{dinh2026}. The lower bound,
\begin{equation}\tag{$L_{n,m}$}
\A_{n,m}(e^\alpha,e^\beta)\ \ge\ \Tr e^{n\alpha+m\beta},
\end{equation}
is open as far as we know. A proof has to use the average over all words: single words in two positive definite letters can have negative trace \cite{hillarjohnson2002}, so $\Tr W\ge\Tr e^{n\alpha+m\beta}$ fails for individual words $W$.

All proofs below go through the following stronger statement for positive semidefinite $A,B$, which compares the average with the single ``evenly spread'' word $(A^{n/m}B)^m$:
\begin{equation}\tag{$E_{n,m}$}
\A_{n,m}(A,B)\ \ge\ \Tr\big(A^{n/m}B\big)^m .
\end{equation}
Here and below $A^s$ is defined by the spectral theorem. By the Araki--Lieb--Thirring inequality, $(E_{n,m})$ implies $(L_{n,m})$ (Lemma~\ref{lem:EL}).

\begin{theorem}\label{thm:main}
Let $n,m\ge0$ and let $\alpha,\beta$ be Hermitian $d\times d$ matrices. The lower bound $(L_{n,m})$ holds in each of the following cases:
\begin{enumerate}
\item[(a)] $\min(n,m)\le2$, any $d$;
\item[(b)] $d\le2$, any $n,m$;
\item[(c)] $d=3$ and $n+m\le7$.
\end{enumerate}
In particular $(L_{n,m})$ holds for all $n+m\le7$ whenever $d\le3$.
\end{theorem}

\begin{theorem}\label{thm:E}
For positive semidefinite $d\times d$ matrices $A,B$ the inequality $(E_{n,m})$ holds if $1\le m\le2$ (any $d$, $n$), if $d\le2$ (any $n$ and $m\ge1$), and if $d=3$, $m=3$, $n\in\{3,4\}$. In particular, for positive semidefinite $3\times3$ matrices,
\[
\tfrac1{20}\Tr S_{6,3}(A,B)\ge\Tr(AB)^3,\qquad \tfrac1{35}\Tr S_{7,3}(A,B)\ge\Tr\big(A^{4/3}B\big)^3 .
\]
\end{theorem}

The cases $m\le2$ and $d\le2$ are elementary: in an eigenbasis of $A$ all weights of the path expansion are nonnegative, and the inequality of arithmetic and geometric means (AM--GM) applies path by path (Section~\ref{sec:paths}). With three letters $B$ in dimension three, only one kind of term can be negative: the ``triangle'' terms with weight $B_{12}B_{23}B_{31}$ or its complex conjugate. An AM--GM pairing, which uses that a $3\times3$ matrix has only one index triangle, removes $B$ and leaves an inequality between the three eigenvalues of $A$ (Section~\ref{sec:three}). Two elementary estimates reduce it to the nonnegativity of one polynomial in two variables on a triangle. For $n=3$ and $n=4$ we certify this in exact rational arithmetic by Bernstein expansion on three boxes (Section~\ref{sec:cert}). The program, which uses only the Python standard library and runs in a few seconds, is provided as an ancillary file.

Theorem~\ref{thm:main} follows from Theorem~\ref{thm:E} because $(L_{n,m})\Leftrightarrow(L_{m,n})$ and because for $d\le3$ and $n+m\le7$ either $\min(n,m)\le2$, or $d\le2$, or $\{n,m\}=\{3,3\}$ or $\{3,4\}$ (Section~\ref{sec:proofs}). In general dimension the first open case is $(n,m)=(3,3)$.

\section{Preliminaries}\label{sec:prelim}

\begin{lemma}[symmetry]\label{lem:sym}
$\A_{n,m}(A,B)=\A_{m,n}(B,A)$. Hence $(L_{n,m})$ for $(\alpha,\beta)$ is $(L_{m,n})$ for $(\beta,\alpha)$.
\end{lemma}
\begin{proof}
Both sides are the average of $\Tr W$ over the same set of words, namely all words with $n$ letters $A$ and $m$ letters $B$. \qed
\end{proof}

\begin{lemma}[from $(E)$ to $(L)$]\label{lem:EL}
For Hermitian $\alpha,\beta$ and $m\ge1$, $\Tr\big(e^{n\alpha/m}e^{\beta}\big)^m\ge\Tr e^{n\alpha+m\beta}$. Consequently $(E_{n,m})$ for $A=e^\alpha$, $B=e^\beta$ implies $(L_{n,m})$.
\end{lemma}
\begin{proof}
For Hermitian $H,K$ the function $r\mapsto\Tr\big(e^{rK/2}e^{rH}e^{rK/2}\big)^{1/r}$ is nondecreasing on $(0,\infty)$; this follows from the Araki--Lieb--Thirring inequality $\Tr(Y^{r/2}X^rY^{r/2})^q\le\Tr(Y^{1/2}XY^{1/2})^{rq}$ for $X,Y\ge0$, $0<r\le1$, $q>0$ \cite{araki1990,liebthirring1976} (see also \cite{carlenlieb2017}). As $r\downarrow0$ it tends to $\Tr e^{H+K}$ by the Lie--Trotter formula \cite{trotter1959}. Take $r=1/m$, $H=n\alpha$, $K=m\beta$ and use cyclicity: $\Tr(e^{n\alpha/m}e^\beta)^m=\Tr(e^{\beta/2}e^{n\alpha/m}e^{\beta/2})^m\ge\Tr e^{n\alpha+m\beta}$. \qed
\end{proof}

\begin{lemma}[average over compositions]\label{lem:comp}
Let $m\ge1$ and let $\C_{n,m}=\{c\in\N_0^m:\ c_1+\dots+c_m=n\}$ be the set of compositions of $n$ into $m$ parts, $|\C_{n,m}|=\binom{n+m-1}{m-1}$. For any square matrices $A,B$ of equal size put $t(c)=\Tr\big(BA^{c_1}BA^{c_2}\cdots BA^{c_m}\big)$. Then
\[
\A_{n,m}(A,B)=\frac1{|\C_{n,m}|}\sum_{c\in\C_{n,m}}t(c).
\]
\end{lemma}
\begin{proof}
A word with $m$ letters $B$ is $W=A^{c_0}BA^{c_1}B\cdots A^{c_{m-1}}BA^{c'_m}$ with $c_0+c_1+\dots+c_{m-1}+c'_m=n$. By cyclicity $\Tr W=t(c_1,\dots,c_{m-1},c_m)$ with $c_m=c_0+c'_m$, and each $c\in\C_{n,m}$ arises from exactly $c_m+1$ words. Since $t$ is invariant under cyclic shifts of $c$ and $\C_{n,m}$ is invariant under them, we may replace the weight $c_m+1$ by its average over the $m$ cyclic shifts, $\frac nm+1$. Hence $\Tr S_{n+m,m}=\frac{n+m}m\sum_ct(c)$, and $\frac{n+m}m\binom{n+m-1}{m-1}=\binom{n+m}m$. \qed
\end{proof}

\section{The path expansion and nonnegative weights}\label{sec:paths}
Let $A=\mathrm{diag}(a_1,\dots,a_d)$ be the diagonal matrix with entries $a_i\ge0$ (we use $0^0=1$) and let $B$ be Hermitian. Expanding $t(c)$ in this basis,
\begin{equation}\label{eq:path}
t(c)=\sum_{i_1,\dots,i_m=1}^d w(i)\prod_{j=1}^m a_{i_{j+1}}^{c_j},\qquad w(i)=B_{i_1i_2}B_{i_2i_3}\cdots B_{i_mi_1},
\end{equation}
with $i_{m+1}:=i_1$. The weight $w(i)$ belongs to the closed index path $i_1\to i_2\to\dots\to i_m\to i_1$. Similarly $\Tr(A^{n/m}B)^m=\sum_iw(i)\prod_ja_{i_j}^{n/m}$. For $x\in[0,\infty)^m$ put
\begin{equation}\label{eq:Delta}
\kappa_{n}(x)=\frac1{|\C_{n,m}|}\sum_{c\in\C_{n,m}}\prod_{j=1}^mx_j^{c_j},\qquad \Delta_{n}(x)=\kappa_n(x)-\Big(\prod_{j=1}^mx_j\Big)^{n/m}.
\end{equation}
Both are symmetric in $x$, because $\C_{n,m}$ is invariant under permutations of the parts. By Lemma~\ref{lem:comp},
\begin{equation}\label{eq:pathdiff}
\A_{n,m}(A,B)-\Tr\big(A^{n/m}B\big)^m=\sum_{i_1,\dots,i_m}w(i)\,\Delta_n\big(a_{i_1},\dots,a_{i_m}\big).
\end{equation}

\begin{lemma}\label{lem:Delta}
$\Delta_n\ge0$ on $[0,\infty)^m$.
\end{lemma}
\begin{proof}
For $x_j>0$ the inequality of arithmetic and geometric means (AM--GM), $\frac1K\sum_{l=1}^Kz_l\ge(\prod_{l=1}^Kz_l)^{1/K}$ for $z_l\ge0$, applied to the $|\C_{n,m}|$ numbers $\prod_jx_j^{c_j}$ gives $\kappa_n(x)\ge\prod_jx_j^{\bar c_j}$ with $\bar c_j=|\C_{n,m}|^{-1}\sum_cc_j=n/m$ (by symmetry the $m$ averages are equal and sum to $n$). Zero arguments follow by continuity. \qed
\end{proof}

\begin{proposition}\label{prop:nonneg}
If all path weights $w(i)$ in an eigenbasis of $A$ are nonnegative, then $(E_{n,m})$ holds. This is the case if $B\ge0$ and $m\le2$ (any $d$), or if $B\ge0$ and $d\le2$ (any $m$).
\end{proposition}
\begin{proof}
The first statement follows from \eqref{eq:pathdiff} and Lemma~\ref{lem:Delta}. For $B\ge0$: if $m=1$, $w(i)=B_{i_1i_1}\ge0$; if $m=2$, $w(i)=B_{i_1i_2}B_{i_2i_1}=|B_{i_1i_2}|^2$. If $d=2$, a closed path passes from $1$ to $2$ as often as from $2$ to $1$, so $w(i)=\prod(\text{diagonal entries})\cdot|B_{12}|^{2k}\ge0$. \qed
\end{proof}

\section{Three letters $B$ in dimension three}\label{sec:three}
From now on $d=3$, $m=3$, and $n\ge1$. We write $\C=\C_{n,3}$, $N=|\C|=\binom{n+2}2$, and $\Delta=\Delta_n$, so that
\[
\Delta(x,y,z)=\frac1N\sum_{c\in\C}x^{c_1}y^{c_2}z^{c_3}-(xyz)^{n/3}.
\]
$\Delta$ is symmetric, continuous on $[0,\infty)^3$ and homogeneous of degree $n$, and $\Delta\ge0$ (Lemma~\ref{lem:Delta}). Moreover $\Delta(t,t,1)>0$ for $0\le t<1$: equality in the AM--GM inequality would require all numbers $t^{c_1+c_2}\in\{1,t,\dots,t^n\}$ to be equal; and $\Delta(0,0,1)=1/N$.

\begin{proposition}[exact identity]\label{prop:id}
Let $A=\mathrm{diag}(a_1,a_2,a_3)$, $a_i\ge0$, and $B$ Hermitian. Then, with $\mathrm{Re}$ denoting the real part,
\begin{multline}\label{eq:id}
\A_{n,3}(A,B)-\Tr\big(A^{n/3}B\big)^3\\
=3\sum_{p\ne q}B_{pp}|B_{pq}|^2\,\Delta(a_p,a_p,a_q)+6\,\mathrm{Re}(B_{12}B_{23}B_{31})\,\Delta(a_1,a_2,a_3).
\end{multline}
\end{proposition}
\begin{proof}
Sort the index triples in \eqref{eq:pathdiff}. If $i_1=i_2=i_3$ then $\Delta(a,a,a)=0$. If the value $p$ occurs twice and $q\ne p$ once, the three triples $(p,p,q)$, $(p,q,p)$, $(q,p,p)$ all have weight $B_{pp}B_{pq}B_{qp}=B_{pp}|B_{pq}|^2$ and, by symmetry, kernel $\Delta(a_p,a_p,a_q)$. If the indices are distinct, the three cyclic orders of $(1,2,3)$ have weight $T:=B_{12}B_{23}B_{31}$, the three others $\overline T$, and all have kernel $\Delta(a_1,a_2,a_3)$. \qed
\end{proof}

For $B\ge0$ all terms of the first sum are nonnegative; only the triangle term can be negative. For $x=(x_1,x_2,x_3)\in[0,\infty)^3$ let $\Delta_{pq}(x)=\Delta(x_p,x_p,x_q)$ and
\begin{equation}\label{eq:G}
G_n(x)=\sum_{k=1}^3\sqrt{\Delta_{ik}(x)\Delta_{jk}(x)}\ -\ \Delta(x_1,x_2,x_3),\qquad \{i,j,k\}=\{1,2,3\}.
\end{equation}
If two entries of $x$ coincide, say $x_1=x_2$, then $G_n(x)=0$: the term $k=3$ equals $\Delta(x_1,x_2,x_3)$, and the other two terms contain $\Delta(x_1,x_1,x_1)=0$.

\begin{lemma}[pairing]\label{lem:pair}
Let $A=\mathrm{diag}(a)$, $a\in[0,\infty)^3$, and $B\ge0$. If $G_n(a)\ge0$, then the right-hand side of \eqref{eq:id} is $\ge0$.
\end{lemma}
\begin{proof}
If some $B_{pp}=0$, the $p$-th row and column of $B$ vanish, $T=0$, and all remaining terms are $\ge0$. Otherwise let $\beta_p=B_{pp}>0$, $\Pi=\beta_1\beta_2\beta_3$ and $u_{pq}=|B_{pq}|/\sqrt{\beta_p\beta_q}\in[0,1]$ (the $2\times2$ principal minors of $B$ are $\ge0$). For $p\ne q$ and $r$ the third index, $B_{pp}|B_{pq}|^2=\Pi\,u_{pq}^2\beta_p/\beta_r$, and $\mathrm{Re}\,T\ge-|T|=-\Pi u_{12}u_{13}u_{23}$. Group the six pairs by their second index $k$; the two pairs $(i,k)$ and $(j,k)$ carry the reciprocal ratios $\beta_i/\beta_j$ and $\beta_j/\beta_i$. By AM--GM for two numbers, $s+t\ge2\sqrt{st}$, and by $u_{ij}\le1$,
\[
u_{ik}^2\Delta_{ik}\frac{\beta_i}{\beta_j}+u_{jk}^2\Delta_{jk}\frac{\beta_j}{\beta_i}\ \ge\ 2u_{ik}u_{jk}\sqrt{\Delta_{ik}\Delta_{jk}}\ \ge\ 2u_{12}u_{13}u_{23}\sqrt{\Delta_{ik}\Delta_{jk}}.
\]
Summing over $k$, the right-hand side of \eqref{eq:id} is $\ge6\,\Pi\,u_{12}u_{13}u_{23}\,G_n(a)\ge0$. \qed
\end{proof}

$G_n$ is symmetric, continuous and homogeneous of degree $n$, and vanishes when two entries coincide. Hence $G_n\ge0$ on $[0,\infty)^3$ follows from $G_n(x,y,1)\ge0$ for $0\le x\le y<1$. For such $(x,y)$ define
\begin{gather*}
\delta_x=\Delta(x,x,1),\qquad\delta_y=\Delta(y,y,1),\qquad\Delta_0=\Delta(x,y,1),\\
Q=\Delta(y,y,x)\Delta(1,1,x)+\Delta(x,x,y)\Delta(1,1,y),\qquad D=\Delta_0^2-\delta_x\delta_y,\\
T'=\Delta_0(\delta_x+\delta_y)+2\delta_x\delta_y,\qquad P_n(x,y)=Q\,T'^2-D^2(\delta_x+\delta_y)^2 .
\end{gather*}
With $(x_1,x_2,x_3)=(x,y,1)$ the square roots in \eqref{eq:G} are $S_3=\sqrt{\delta_x\delta_y}$, $S_1=\sqrt{\Delta(y,y,x)\Delta(1,1,x)}$, $S_2=\sqrt{\Delta(x,x,y)\Delta(1,1,y)}$, and $Q=S_1^2+S_2^2$.

\begin{lemma}[removing the square roots]\label{lem:sqrt}
If $P_n(x,y)\ge0$ for all $0\le x\le y<1$, then $G_n\ge0$ on $[0,\infty)^3$.
\end{lemma}
\begin{proof}
Since $S_1+S_2\ge\sqrt Q$, we have $G_n(x,y,1)\ge\sqrt Q-(\Delta_0-S_3)$. If $D\le0$, then $S_3\ge\Delta_0$. If $D>0$, then $\Delta_0>0$, and $\delta_y>0$ because $y<1$. The inequality between geometric and harmonic mean, $S_3=\sqrt{\delta_x\delta_y}\ge2\delta_x\delta_y/(\delta_x+\delta_y)$ (which follows from $\delta_x+\delta_y\ge2\sqrt{\delta_x\delta_y}$), gives
\[
\Delta_0-S_3=\frac{D}{\Delta_0+S_3}\le\frac{D}{\Delta_0+2\delta_x\delta_y/(\delta_x+\delta_y)}=\frac{D(\delta_x+\delta_y)}{T'} .
\]
So $G_n\ge0$ if $\sqrt Q\,T'\ge D(\delta_x+\delta_y)$; as both sides are $\ge0$, this is $P_n\ge0$. \qed
\end{proof}

The cruder bound $S_3\ge0$ in place of the harmonic mean would lead to the condition $Q\Delta_0^2\ge D^2$, which fails on the edge $x=0$ for $y\gtrsim0.45$, both for $n=3$ and for $n=4$.

\section{The certificate for $n=3$ and $n=4$}\label{sec:cert}
Put $k=1$ if $n=3$ and $k=3$ if $n=4$, and substitute $x=X^k$, $y=Y^k$ with $0\le X\le Y\le1$. Then $(xyz)^{n/3}=(XYZ)^{kn/3}$ is a monomial, every value of $\Delta$ occurring above is a polynomial in $(X,Y)$ with rational coefficients, and so is $\Pc_n(X,Y):=P_n(X^k,Y^k)$. It vanishes on the edges $X=Y$ (where $D=Q=0$) and $Y=1$ (where $\delta_y=0$ and $Q=\Delta_0^2$) of the triangle
\[
\Theta=\{(X,Y):\ 0\le X\le Y\le1\}.
\]

\begin{lemma}\label{lem:cert}
For $n\in\{3,4\}$: \textup{(a)} $\Pc_n=(Y-X)^2(1-Y)^2R_n$ with $R_n\in\Q[X,Y]$; \textup{(b)} $R_n\ge0$ on $\Theta$.
\end{lemma}

Table~\ref{tab:data} lists the data. Part (a) is an exact polynomial division. For (b) we use the \emph{Bernstein principle}: if $p(s,t)=\sum_{i\le\mu,j\le\nu}b_{ij}\binom\mu is^i(1-s)^{\mu-i}\binom\nu jt^j(1-t)^{\nu-j}$, then $p\ge\min_{ij}b_{ij}$ on $[0,1]^2$, since the basis functions are nonnegative and sum to $1$. The coefficients $b_{ij}$ on a box are obtained by the affine change of variables onto $[0,1]^2$ followed by $b_l=\sum_{i\le l}\binom li\binom\mu i^{-1}a_i$ in each variable.

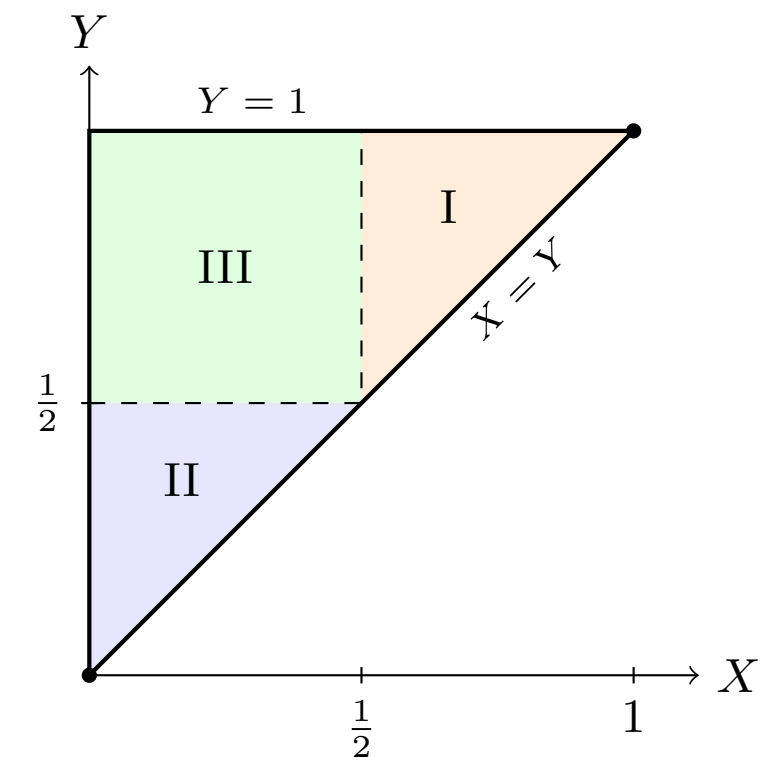
\begin{figure}
\centering
\begin{tikzpicture}[scale=3.6]
\fill[blue!10] (0,0)--(0.5,0.5)--(0,0.5)--cycle;
\fill[green!12] (0,0.5)--(0.5,0.5)--(0.5,1)--(0,1)--cycle;
\fill[orange!15] (0.5,0.5)--(1,1)--(0.5,1)--cycle;
\draw[thick] (0,0)--(1,1)--(0,1)--cycle;
\draw[dashed] (0,0.5)--(0.5,0.5)--(0.5,1);
\draw[->] (0,0)--(1.12,0) node[right]{$X$};
\draw[->] (0,0)--(0,1.12) node[above]{$Y$};
\draw (0.5,0.015)--(0.5,-0.015) node[below]{$\frac12$};
\draw (1,0.015)--(1,-0.015) node[below]{$1$};
\draw (0.015,0.5)--(-0.015,0.5) node[left]{$\frac12$};
\node at (0.17,0.36) {II};
\node at (0.25,0.75) {III};
\node at (0.66,0.86) {I};
\fill (0,0) circle (0.014); \fill (1,1) circle (0.014);
\node[rotate=45,below] at (0.75,0.75) {\scriptsize $X=Y$};
\node[above] at (0.3,1) {\scriptsize $Y=1$};
\end{tikzpicture}
\caption{The triangle $\Theta$ and the three charts of Lemma~\ref{lem:cert}. $\Pc_n$ vanishes on the edges $X=Y$ and $Y=1$ (the factor $(Y-X)^2(1-Y)^2$); $R_n$ vanishes only at the two marked corners, where charts II and I divide out the zero.}
\label{fig:charts}
\end{figure}

$R_n$ vanishes at the corners $(0,0)$ and $(1,1)$ of $\Theta$, where the corresponding vertex coefficient is $0$. We therefore use three charts (Fig.~\ref{fig:charts}):
\begin{description}
\item[I ($X\ge\frac12$):] $a=Y-X=st$, $b=1-Y=s(1-t)$ with $(s,t)\in[0,\frac12]\times[0,1]$; here $s=1-X$ and $t=a/s$. Written in $(a,b)$, every monomial of $R_n$ has degree $\ge\omega_1$, so $R_{n,\rm I}:=R_n/s^{\omega_1}$ is a polynomial in $(s,t)$.
\item[II ($Y\le\frac12$):] $X=rY$ with $(r,Y)\in[0,1]\times[0,\frac12]$. Every monomial of $R_n$ has degree $\ge\omega_0$, so $R_{n,\rm II}:=R_n(rY,Y)/Y^{\omega_0}$ is a polynomial.
\item[III ($X\le\frac12\le Y$):] $R_n$ itself on $[0,\frac12]\times[\frac12,1]$.
\end{description}
If $Y\le\frac12$ we are in II; if $X\le\frac12\le Y$, in III; if $X\ge\frac12$, in I. All Bernstein coefficients of the three chart polynomials on their boxes are positive (Table~\ref{tab:data}), so $R_n\ge0$ on $\Theta$. This proves Lemma~\ref{lem:cert}.

\begin{table}
\caption{Data of the certificate (computed by the ancillary program in exact rational arithmetic). $\omega_0,\omega_1$: orders of the zeros of $R_n$ at $(0,0)$ and $(1,1)$.}
\label{tab:data}
\centering\small
\begin{tabular}{@{}p{3.6cm}p{3.3cm}p{2.6cm}@{}}
\toprule
 & $n=3$ (word sum $S_{6,3}$) & $n=4$ ($S_{7,3}$)\\
\midrule
substitution & $x=X$ & $x=X^3$\\
degree of $\Pc_n$ / of $R_n$ & $18$ / $14$ & $72$ / $68$\\
nonzero coefficients of $R_n$ & $116$ & $2283$\\
$\omega_0$, $\omega_1$ & $1$, $8$ & $10$, $8$\\
min.\ coefficient, chart I & $\approx1.4\cdot10^{-2}$ & $\approx3.3\cdot10^{-3}$\\
min.\ coefficient, chart II & $729/3200000$ & $32/759375$\\
min.\ coefficient, chart III & $\approx5.4\cdot10^{-5}$ & $\approx2.4\cdot10^{-7}$\\
\bottomrule
\end{tabular}
\end{table}

\paragraph{Reproducibility.} The ancillary file \texttt{bmv\_certificate.py} (Python~3, standard library only) builds $\Pc_n$ from the definitions above, performs the division of Lemma~\ref{lem:cert}(a) (a nonzero remainder stops the program), forms the three chart polynomials and prints their smallest Bernstein coefficients as exact fractions; \texttt{python3 bmv\_certificate.py 3} and \texttt{\dots\ 4} run in less than a second and about three seconds, respectively. The ancillary file \texttt{sympy\_\allowbreak crosscheck.py} repeats the whole construction independently with the Python computer-algebra library SymPy (polynomial division, chart substitutions, Bernstein coefficients); for $n=3$ and $n=4$ it reproduces $R_n$ coefficient by coefficient and the minima of Table~\ref{tab:data}. The ancillary file \texttt{numerical\_checks.py} contains floating-point checks of \eqref{eq:id} and of Theorem~\ref{thm:E}; they are not part of the proof.

\section{Proofs of the theorems}\label{sec:proofs}
\begin{proof}[of Theorem~\ref{thm:E}]
Both sides of $(E_{n,m})$ are invariant under $(A,B)\mapsto(UAU^*,UBU^*)$ with unitary $U$, so we may take $A$ diagonal. The cases $m\le2$ and $d\le2$ are Proposition~\ref{prop:nonneg}. For $d=3$, $m=3$, $n\in\{3,4\}$: Lemma~\ref{lem:cert} gives $\Pc_n\ge0$ on $\Theta$, that is, $P_n(x,y)\ge0$ for $0\le x\le y\le1$; Lemma~\ref{lem:sqrt} gives $G_n\ge0$; Lemma~\ref{lem:pair} and Proposition~\ref{prop:id} give $(E_{n,3})$. Finally $\binom63=20$, $\binom73=35$. \qed
\end{proof}

\begin{proof}[of Theorem~\ref{thm:main}]
If $d=1$ or $\min(n,m)=0$, $(L_{n,m})$ is an equality. Let $A=e^\alpha$, $B=e^\beta$ and $n,m\ge1$. If $m\le2$, or $d=2$, Theorem~\ref{thm:E} and Lemma~\ref{lem:EL} give $(L_{n,m})$; if $n\le2$, apply this to $(L_{m,n})$ and use Lemma~\ref{lem:sym}. This proves (a) and (b). For (c) it remains to treat $d=3$, $n+m\le7$, $\min(n,m)\ge3$, that is, $(n,m)\in\{(3,3),(4,3),(3,4)\}$. The first two follow from Theorem~\ref{thm:E} and Lemma~\ref{lem:EL}, the third from the second by Lemma~\ref{lem:sym}. \qed
\end{proof}

\section{Remarks}\label{sec:remarks}
\textit{(i) Where the argument loses.} The proof of Lemma~\ref{lem:pair} uses $\mathrm{Re}\,T\ge-|T|$, the AM--GM pairing and $u_{ij}\le1$; these are equalities for commuting $A,B$. The two estimates of Lemma~\ref{lem:sqrt} do not involve $B$; their loss is multiplied by $u_{12}u_{13}u_{23}$, which vanishes for commuting matrices. The margins of the certificate are small but genuine: for $n=4$ the smallest coefficient in chart III is the vertex value $R_4(0,\frac12)$.

\textit{(ii) Other numbers of letters.} Proposition~\ref{prop:id}, Lemma~\ref{lem:pair} and Lemma~\ref{lem:sqrt} hold for every $n$. Thus $(E_{n,3})$ for $d=3$, and with it $(L_{n,3})$ and $(L_{3,n})$, follows from $P_n\ge0$ on $0\le x\le y<1$. We have not treated $n\ge5$ here.

\textit{(iii) Higher dimension.} Lemma~\ref{lem:pair} uses that in dimension three every edge $\{p,q\}$ belongs to exactly one index triangle. For $d\ge4$ several triangles share edges and the pairing does not carry over directly. The cases $(3,3)$ and $(4,3)$ remain open for $d\ge4$; for $d=3$ the first cases not covered here are $(4,4)$ and $(5,3)$, $(3,5)$.

\begin{acknowledgements}
Artificial intelligence (AI) consulting: H\"agele Consulting. The proofs were checked independently by a Claude Opus 5.5 assistant (Anthropic) and by a Gemini Argon agent (Google). An independent verification of the proofs by a human referee is still pending.
\end{acknowledgements}

\section*{Declarations}
\textbf{Use of generative AI.} The derivations, the computer certificate and a draft of the text were developed with the assistance of the large language model Claude Opus 5.5 (Anthropic), guided by the author. The author has reviewed the manuscript and takes full responsibility for its content.\\
\textbf{Funding.} No funding was received for this work.\\
\textbf{Competing interests.} The author declares no competing interests.\\
\textbf{Data and code availability.} The programs \texttt{bmv\_\allowbreak certificate.py}, \texttt{sympy\_\allowbreak crosscheck.py} and \texttt{numerical\_\allowbreak checks.py} are available as ancillary files of the arXiv version of this article.

\bibliographystyle{spmpsci}
\bibliography{bmv_lower}

\begin{thebibliography}{10}
\providecommand{\url}[1]{{#1}}
\providecommand{\urlprefix}{URL }
\expandafter\ifx\csname urlstyle\endcsname\relax
  \providecommand{\doi}[1]{DOI~\discretionary{}{}{}#1}\else
  \providecommand{\doi}{DOI~\discretionary{}{}{}\begingroup
  \urlstyle{rm}\Url}\fi

\bibitem{araki1990}
Araki, H.: On an inequality of {Lieb} and {Thirring}.
\newblock Lett. Math. Phys. \textbf{19}, 167--170 (1990)

\bibitem{bmv1975}
Bessis, D., Moussa, P., Villani, M.: Monotonic converging variational
  approximations to the functional integrals in quantum statistical mechanics.
\newblock J. Math. Phys. \textbf{16}(11), 2318--2325 (1975)

\bibitem{carlenlieb2017}
Carlen, E.A., Lieb, E.H.: Some trace inequalities for exponential and
  logarithmic functions.
\newblock arXiv:1709.05450 (2017)

\bibitem{cha2026}
Cha, H., Lee, J.: One-parameter counterexamples to the refined
  {Bessis--Moussa--Villani} conjecture.
\newblock arXiv:2603.19927v6 (2026)

\bibitem{dinh2026}
Dinh, T.H.: On the failure of the upper bound in the refined {BMV} conjecture
  and a pinching correction.
\newblock arXiv:2605.17782 (2026)

\bibitem{golden1965}
Golden, S.: Lower bounds for the {Helmholtz} function.
\newblock Phys. Rev. \textbf{137}, B1127--B1128 (1965)

\bibitem{haegele2007}
H{\"a}gele, D.: Proof of the cases $p \le 7$ of the {Lieb--Seiringer}
  formulation of the {Bessis--Moussa--Villani} conjecture.
\newblock J. Stat. Phys. \textbf{127}, 1167--1171 (2007)

\bibitem{hillarjohnson2002}
Hillar, C.J., Johnson, C.R.: Eigenvalues of words in two positive definite
  letters.
\newblock SIAM J. Matrix Anal. Appl. \textbf{23}, 916--928 (2002).
\newblock ArXiv:math/0511411

\bibitem{klepschweighofer2008}
Klep, I., Schweighofer, M.: Sums of {H}ermitian squares and the {BMV}
  conjecture.
\newblock J. Stat. Phys. \textbf{133}, 739 (2008)

\bibitem{liebseiringer2004}
Lieb, E.H., Seiringer, R.: Equivalent forms of the {Bessis--Moussa--Villani}
  conjecture.
\newblock J. Stat. Phys. \textbf{115}, 185--190 (2004)

\bibitem{liebthirring1976}
Lieb, E.H., Thirring, W.: Inequalities for the moments of the eigenvalues of
  the {Schr\"odinger} {Hamiltonian} and their relation to {Sobolev}
  inequalities.
\newblock In: E.H. Lieb, B.~Simon, A.S. Wightman (eds.) Studies in Mathematical
  Physics, pp. 269--303. Princeton University Press (1976)

\bibitem{oqp40}
{Open Quantum Problems}: Problem 40: Refinement of the
  {Bessis--Moussa--Villani} conjecture.
\newblock
  \url{https://oqp.iqoqi.oeaw.ac.at/refinement-of-the-bessis-moussa-villani-conjecture}
  (accessed 2026)

\bibitem{stahl2013}
Stahl, H.R.: Proof of the {BMV} conjecture.
\newblock Acta Math. \textbf{211}, 255--290 (2013)

\bibitem{thompson1965}
Thompson, C.J.: Inequality with applications in statistical mechanics.
\newblock J. Math. Phys. \textbf{6}, 1812--1813 (1965)

\bibitem{trotter1959}
Trotter, H.F.: On the product of semi-groups of operators.
\newblock Proc. Amer. Math. Soc. \textbf{10}, 545--551 (1959)

\end{thebibliography}

\end{document}